\documentstyle[amssymb,amsfonts]{amsart}

\def\bi{\begin{itemize}}
\def\bs{\begin{split}}
\def\es{\end{split}}
\def\ba{\begin{align}}
\def\bas{\begin{align*}}
\def\ea{\end{align}}
\def\eas{\end{align*}}
\def\R{{\hbox{\bf R}}}

\def\R{{\hbox{\bf R}}}
\def\C{{\hbox{\bf C}}}
\def\sgn{{\hbox{sgn}}}

\def\eps{\varepsilon}

\newenvironment{proof}{\noindent {\bf Proof} }{\endprf\par}
\def \endprf{\hfill  {\vrule height6pt width6pt depth0pt}\medskip}
\def\emph#1{{\it #1}}
\def\textbf#1{{\bf #1}}

\theoremstyle{plain}
  \newtheorem{theorem}[subsection]{Theorem}
  
  \newtheorem{proposition}[subsection]{Proposition}
  \newtheorem{lemma}[subsection]{Lemma}
  \newtheorem{corollary}[subsection]{Corollary}

\theoremstyle{remark}

\theoremstyle{definition}
  \newtheorem{definition}[subsection]{Definition}

\include{psfig}

\begin{document}

\title[Well-posedness of Benjamin-Ono equation]{Global well-posedness of the Benjamin-Ono equation in $H^1(\R)$}
\author{Terence Tao}
\address{Department of Mathematics, UCLA, Los Angeles CA 90095-1555}
\email{ tao@@math.ucla.edu}
\subjclass{35J10}

\vspace{-0.3in}
\begin{abstract}
We show that the Benjamin-Ono equation is globally well-posed in $H^s(\R)$ for $s \geq 1$.  This is despite the presence of the derivative in the non-linearity, which causes the solution map to not be uniformly continuous in $H^s$ for any $s$ \cite{kt-2}.  The main new ingredient is to perform a global gauge transformation which almost entirely eliminates this derivative.  
\end{abstract}

\maketitle

\section{Introduction}

In this paper we study the Cauchy problem for the Benjamin-Ono equation
\begin{equation}\label{eq:benjamin}
 u_t + H u_{xx} = u u_x; \quad u(0,x) = u_0(x),
\end{equation}
where $u: \R \times \R \to \R$ is a real-valued function, $H$ is the spatial Hilbert transform
$$ Hu(x) := \frac{1}{\pi} p.v. \int_\R \frac{u(y)}{x-y}\ dy,$$
and $u_0$ lives in the inhomogeneous Sobolev space $H^s_x$ for some $s \in \R$.  Recall that these spaces can be defined via the spatial Fourier transform
$$ \hat f(\xi) := \int_{\R} e^{-2\pi i x \xi} f(x)\ dx$$
by the formula
$$ \| f \|_{H^s_x} := (\int_{\R} \langle \xi \rangle^{2s} |\hat f(\xi)|^2\ d\xi)^{1/2},$$
where $\langle \xi \rangle := (1 + |\xi|^2)^{1/2}$.

This equation is a model for one-dimensional waves in deep water \cite{bo}, and is completely integrable \cite{af}.  The Cauchy problem for this equation has been extensively studied \cite{gv}, \cite{gv-2}, \cite{iorio}, \cite{kt}, \cite{kt-2}, \cite{ponce}, \cite{saut}, \cite{tom}.   
It is known that this equation is globally well-posed in the Sobolev space $H^s_x$ if $s$ is sufficiently large; for instance it is known that if $u_0 \in H^3_x$ (for instance), then there exists a unique (classical) global solution $u$ to \eqref{eq:benjamin} which lies in the space $C^0_t H^3_x$; we then define the solution operator $S(t): H^3_x \to H^3_x$.  Furthermore, it is known that for every fixed $t$, the map $S(t): H^3_x \to H^3_x$ is continuous (see e.g. \cite{saut}), although it is not analytic or even uniformly continuous on bounded sets (see \cite{kt-2}, \cite{mst}).

This paper is concerned with the question of whether $S$ can be extended to rougher classes of initial data.  The above $H^3_x$ result has already been improved by several authors; for instance, Iorio \cite{iorio} obtained local well-posedness in $H^s$ for $s > 3/2$, which was then improved by Ponce \cite{ponce} to well-posedness in $s \geq 3/2$.  More recently, Koch and Tzvetkov \cite{kt} have improved this to $s > 5/4$.

The Hilbert transform behaves in some ways like the constant $\pm i$; for instance, we have that $H$ is anti-self-adjoint with $H^2 = -1$.  A more precise relationship is given by the well-known formula
\begin{equation}\label{eq:signum}
 \widehat{Hf}(\xi) = - i \sgn(\xi) \hat f(\xi)
\end{equation}
which can easily be verified by contour integration.  Thus, the Benjamin-Ono equation \eqref{eq:benjamin} is heuristically of the form
$$ u_t \pm i u_{xx} = u u_x.$$
We remark that quadratic non-linear Schr\"odinger equations of the form
\begin{equation}\label{eq:nls}
 u_t \pm i u_{xx} = |u|^2
\end{equation}
have been studied in \cite{kpv}, and local well-posedness was established for such equations in $H^s(\R)$ for\footnote{It is not known whether this value of $s$ is sharp, however the bilinear estimate used in \cite{kpv} to establish this result certainly fails \cite{kenji} for $s < -1/4$, and indeed can be used to show that the solution map is not analytically well-posed.  See also \cite{taoka}.} $s > -1/4$.  Note that it is not \emph{a priori} obvious that one can make sense (as a distribution) of non-linearities such as $|u|^2$ for solutions in negative Sobolev spaces, but it turns out that the Schr\"odinger equation on the line has sufficient smoothing properties to make this possible.

There is thus a substantial gap - more than a derivative - in regularity between the previous local well-posedness theory for the Benjamin-Ono equation and for the quadratic NLS equation.  We believe that this gap is mostly artificial, in that the well-posedness theory for the Benjamin-Ono equation can be lowered much further, to become closer to that for the quadratic NLS equation.  As a first step in this direction we give

\begin{theorem}\label{thm:lwp}  Let $R > 0$ and $s_0 \geq 1$.  Then there exists a time $T = T(R, s_0) > 0$ such that for every time $-T \leq t \leq T$, the solution operator $S(t): H^3_x \to H^3_x$ can be continuously and uniquely extended to the ball
$$ B(0,R) := \{ u_0 \in H^1_x: \| u_0 \|_{H^1_x} \leq R \},$$
endowed with the $H^1$ topology, and $S(t) u_0$ solves \eqref{eq:benjamin} in the sense of distributions.  In addition, if $u_0 \in B(0,R) \cap H^{s_0}_x$, then $S(t) u_0 \in C^0_{[-T,T]} H^{s_0}_x$, and indeed
\begin{equation}\label{eq:regularity}
 \| S(t) u_0 \|_{C^0_{[-T,T]} H^{s_0}_x} \leq C(s_0,R) \| u_0 \|_{H^{s_0}_x}.
\end{equation}
Here of course we give $C^0_{[-T,T]}$ the sup norm topology:
$$ \| S(t) u_0 \|_{C^0_{[-T,T]} H^{s_0}_x} := \sup_{t \in [-T,T]} \| S(t) u_0 \|_{H^{s_0}_x}.$$
Indeed, the solution map is continuous from $H^{s_0}_x \cap B(0,R)$, with the $H^{s_0}_x$ topology, to $C^0_{[-T,T]} H^{s_0}_x$.

Finally, the solution map is Lipschitz as a map from $B(0,R)$, with the $L^2_x$ topology(!), to $C^0_{[-T,T]} L^2_x$, or in other words that
\begin{equation}\label{lip}
\| S(t) u_0 - S(t) \tilde u_0 \|_{L^2_x} \leq C(s_0,R) \| u_0 - \tilde u_0 \|_{L^2_x}
\end{equation}
for all $t \in [-T,T]$ and $u_0, \tilde u_0 \in B(0,R)$.
\end{theorem}

This improves\footnote{Since the submission of this paper, the author has learned (C. Kenig, personal communication)
that a similar improvement by different methods has been achieved by C. Kenig and K. Koenig, for $s > 9/8$.}
on the earlier result of Koch and Tzvetkov \cite{kt}, which obtained local well-posedness in $H^s_x$ for $s > 5/4$.  Also, in \cite{kt-2} it is remarked that the solution map is not uniformly continuous from $H^s_x$ to $H^s_x$ for any $s$, because of the derivative in the non-linearity and the relatively weak smoothing effects of the linear part of the equation, and so a direct iteration map cannot work; note however that the construction in \cite{kt-2} does not prohibit the solution map being uniformly continuous or Lipschitz in a weaker topology such as $L^2_x$.  We circumvent this lack of uniform continuity by applying a gauge transform (a variant of the classical Cole-Hopf transformation) to effectively remove the derivative (or at least the worst terms which involve\footnote{Specifically, we can eliminate the terms involving an interaction of very low and very high frequencies, where the derivative falls on the very high frequency; this interaction was specifically identified as the culprit to lack of uniform continuity in \cite{kt-2}.} the derivative) from the non-linearity.  This gauge transform may well be related to the complete integrability of the Benjamin-Ono equation (see e.g. \cite{af}), although we do not know of a direct connection between our methods and the theory of integrable systems.

Our methods used are fairly elementary, relying mostly on the algebraic manipulations of the gauge transform, 
as well as standard Strichartz estimates for the Schr\"odinger equation.  
Littlewood-Paley projections play a small role (except in obtaining the $H^1$ continuity estimate), 
but mostly we need them to separate the action of the Hilbert transform on the positive and negative frequencies.  
The Lipschitz bound \eqref{lip} in particular is very simple and relies on the energy method combined with 
the $L^4_t L^\infty_x$ Strichartz estimate.

For smooth (e.g. $H^3_x$) solutions, it is easy to check that the $L^2_x$ norm $\int u^2$ is an invariant of the flow \eqref{eq:benjamin}:
\begin{equation}\label{l2-conserv}
 \int_\R u(t,x)^2\ dx = \int_\R u_0(x)^2\ dx.
\end{equation}
Somewhat less obvious conserved quantities include the Hamiltonian
\begin{equation}\label{hamil}
\int_\R u H u_x - \frac{1}{3} u^3\ dx,
\end{equation}
as well as the $H^1$-type quantity
\begin{equation}\label{h1}
 \int_\R u_x^2 - \frac{3}{4} u^2 H u_x - \frac{1}{8} u^4\ dx;
\end{equation}
in fact the Benjamin-Ono equation is completely integrable (see e.g. \cite{af}) and has an infinite number of such conserved quantities.  From the above conservation laws and the Gagliardo-Nirenberg inequality it is easy to verify the \emph{a priori} bound
\begin{equation}\label{u-h1}
 \| u(t) \|_{H^1_x} \leq C( \| u_0 \|_{H^1_x} )
\end{equation}
for all $t \in \R$.  \emph{A priori} this bound can only be proven for smooth solutions, but by the continuity properties obtained in Theorem \ref{thm:lwp} we see that they in fact hold for $H^1_x$ solutions as long as the solution exists.  From this and iterating Theorem \ref{thm:lwp}, we thus obtain

\begin{corollary}\label{thm:gwp}  The equation \eqref{eq:benjamin} is globally well-posed in $H^s_x$ for all $s \geq 1$.
\end{corollary}

Interestingly, a global existence result has recently been obtained for the viscous Benjamin-Ono equation in \cite{f}; it thus seems feasible to show that the solutions to the viscous Benjamin-Ono equation converges to the inviscid Benjamin-Ono equation in the zero viscosity limit, although we do not pursue this question here.

The author is indebted to Nikolay Tzvetkov and Jim Colliander for very helpful discussions on this problem, and to Felipe Linares for some help with the references.  The author also thanks Martin Hadac for pointing out an error in the exposition in a previous version of the paper.

\section{Notation}\label{sec:notation}

Fix $s_0 \geq 1$.  We use the notation $A \lesssim B$ or $A = O(B)$ to denote the estimate that $|A| \leq C_{s_0} B$, where $C_{s_0}$ is a quantity depending on $s_0$ but not on $\eps$.  If $X$ is a Banach space, we use $O_X(B)$ to denote any element in $X$ with norm $O(B)$.  We use $\langle x\rangle$ to denote the quantity $\langle x \rangle := (1 + |x|^2)^{1/2}$.

We observe the Riesz decomposition
$$ 1 = P_- + P_+$$
where $P_\pm$ are the Fourier projections to $\pm [0,\infty)$; from \eqref{eq:signum} we observe that
\begin{equation}\label{iH}
iH = P_- - P_+.
\end{equation}
Let $\psi$ be a bump function adapted to $[-2,2]$ and equal to 1 on $[-1,1]$.  We define the Littlewood-Paley operators $P_{k}$ and $P_{\leq k} = P_{<k+1}$ for $k \geq 0$ by defining 
$$ \widehat{P_{\leq k} f}(\xi) := \psi(\xi/2^k) \hat f(\xi)$$
for all $k \geq 0$, and $P_k := P_{\leq k} - P_{\leq k-1}$ (with the convention $P_{\leq -1} = 0$).  Note that all the operators $P_k$, $P_{\leq k}$ are bounded on all translation-invariant Banach spaces, thanks to Minkowski's inequality.  We define $P_{>k} := P_{\geq k-1} := 1 - P_{\leq k}$.

From Plancherel we have the bound
\begin{equation}\label{eq:planch}
 \| f \|_{H^s_x} \sim (\sum_{k=0}^\infty \| P_k f\|_{H^s_x}^2)^{1/2}
\sim (\sum_{k=0}^\infty 2^{ks} \| P_k f\|_{L^2_x}^2)^{1/2}
\end{equation}
for any $s \in \R$.

We define $P_{lo} := P_0$, $P_{hi} := P_{> 0}$, and $P_{\pm hi} := P_{\pm} P_{hi}$; thus
$1 = P_{lo} + P_{+hi} + P_{-hi}$.   Observe that if $F$ is real, then $P_{lo} F$ is also real; this will be important when we compute expressions such as $e^{i P_{lo} F}$.  We also define $P_{LO} := P_{\leq 2}$, $P_{HI} := P_{>2}$, and $P_{\pm HI} := P_\pm P_{HI}$, so $1 = P_{LO} + P_{+HI} + P_{-HI}$.

\section{The gauge transformation, and a priori Strichartz estimates}

We now begin the proof of Theorem \ref{thm:lwp}.  Fix $s_0$.  We first observe from the scale invariance 
$$ u(t,x) \mapsto \frac{1}{\lambda} u(\frac{t}{\lambda^2}, \frac{x}{\lambda})$$
of the equation \eqref{eq:benjamin}, that we may rescale the $H^1$ norm $R$ to be small.  In particular we may assume that $R \ll \eps^2$, and we will now obtain local well-posedness on the time\footnote{In principle, this scaling should yield a time of existence $T \sim \| u \|_{H^1}^{-4/3}$ when the data is large, since the rescaling shrinks the $\dot H^1$ norm by a factor of $O(\lambda^{-3/2})$.  However our argument does not quite achieve this because the $L^2_x$ component of the $H^1_x$ norm does not scale so well, by a factor of only $O(\lambda^{-1/2})$.  It is quite likely one can modify this argument to fix this problem, although for our purposes this is unnecessary because the $H^1$ conservation law yields a global existence result anyway.} interval $[-1,1]$ (i.e. we set $T := 1$).

From the Strichartz theory of the one-dimensional Schr\"odinger equation (see e.g. \cite{tao:keel}) we have the global Strichartz estimate
$$ \| u \|_{L^4_t L^\infty_x} + \| u \|_{L^\infty_t L^2_x} \lesssim \| u(0) \|_{L^2_x} + \| u_t \pm i u_{xx} \|_{L^1_t L^2_x}$$
for any spacetime test function $u(t,x)$ and either choice of sign $\pm$.  Applying this to $P_\pm u$ and summing, we obtain in particular that
$$ \| u \|_{L^4_t L^\infty_x} + \| u \|_{L^\infty_t L^2_x} \lesssim \| u(0) \|_{L^2_x} + \| u_t + H u_{xx} \|_{L^1_t L^2_x},$$
and similarly for $P_- u$ and $P_+ u$.  Differentiating this, we obtain
\begin{equation}\label{strichartz}
\| u \|_{S^k} \lesssim \| u(0) \|_{H^k_x} + \| u_t + H u_{xx} \|_{L^1_t H^k_x}
\end{equation}
for any integer $k \geq 0$, where $S^k$ denotes the Strichartz norm
$$ \| u \|_{S^k} := \| u \|_{L^4_t C^k_x} + \| u \|_{L^\infty_t H^k_x}.$$
In light of this, it is reasonable to expect that $H^1_x$ solutions of the Benjamin-Ono equation are bounded in $S^1$.  We shall now show that this is in fact the case, at least for smooth solutions.

\begin{theorem}\label{apriori}  Let $u$ be an $H^3_x$ solution to \eqref{eq:benjamin} with $u(0) = O_{H^1}(\eps^2)$.  Then we have
\begin{equation}\label{u-small}
 u = O_{S^1([-1,1] \times \R)}(\eps^2).
\end{equation}
\end{theorem}

\emph{Remark.}  From the classical global well-posedness theory in $H^3_x$ (\cite{saut}; see also \cite{iorio}, \cite{ponce}, \cite{kt}) we then know that as soon as $u_0$ is in $H^3_x$, then the solution $u$ is in $C^0_t H^{3}_x$ and exists globally in time.  Note however that our bounds will be independent of the $H^3_x$ norm of $u_0$, and once we obtain some continuity estimates we will be able to pass to the limit and obtain a version of this theorem for general $H^1_x$ solutions.

\emph{Remark.}  The $L^\infty_t H^1_x$ component of the Strichartz norm $S^1$ can of course be controlled using the conserved quantities of the Benjamin-Ono equation, and in particular \eqref{h1}.  However, we will refrain from using this conservation law here as it does not adapt well to frequency envelopes, which we will introduce in later sections.  We will however take advantage of the $L^2$ conservation law \eqref{l2-conserv} in order to easily deal with the low frequency portions of the solution (although it is possible to treat those without recourse to any conservation law).

\begin{proof}  Henceforth we restrict all spacetime norms to the slab $[-1,1] \times \R$.
A standard iteration method will not work here, because the linear part of the Benjamin-Ono equation does not have enough smoothing to compensate for the derivative in the non-linearity.  To resolve this, we will gauge transform the equation \eqref{eq:benjamin} into a more manageable form, where the derivative term has been moved onto a low frequency factor (cf. the gauge transformations for wave maps in \cite{tao:wavemap}, \cite{tao:wavemap2}, \cite{kr}, \cite{tataru}).  

The first step is to take an antiderivative of \eqref{eq:benjamin}, introducing a spatial primitive $F(t,x)$ of $u(t,x)$.  We first construct $F(t,0)$ on the time axis $x=0$ by solving the ODE
$$ \partial_t F(t,0) + \frac{1}{2} H u_x(t,0) = \frac{1}{4} u(t,0)^2; \quad F(0,0) := 0$$
and then constructing $F(t,x)$ for general $x$ by the ODE
\begin{equation}\label{eq:F-u}
 \partial_x F(t,x) = \frac{1}{2} u(t,x).
\end{equation}
Note that $F(t,x)$ is necessarily real-valued.  From \eqref{eq:benjamin} we see that
$$ \partial_x ( F_t + H F_{xx} ) = \partial_x( F_x^2 )$$
while from our construction on the time axis we know that
$$ F_t(t,0) + H F_{xx}(t,0) = F_x^2(t,0).$$
Thus we have the equation 
\begin{equation}\label{eq:F}
 F_t + H F_{xx} = F_x^2.
\end{equation}
holding globally in spacetime.

{\it Remark.} If we replaced $H$ by $-i$ (as one could do if $F$ were somehow restricted to positive frequencies) then one could solve this equation $F_t - iF_{xx} = F_x^2$ explicitly using the Cole-Hopf transformation\footnote{Actually, this is a complexified version of the classical Cole-Hopf transformation \cite{cole}, \cite{hopf}, used to solve the inviscid Burger's equation.  A similar transformation was also used by Nirenberg (see e.g. \cite{km}) to solve the scalar non-linear wave equation $F_{tt} - \Delta F = |F_t|^2 - |\nabla F|^2$.  We thank Jim Colliander for pointing out the relevance of the Cole-Hopf transformation to this problem.} $w := e^{-iF}$, since one can then easily verify that $w$ solves the free Schr\"odinger equation $w_t - iw_{xx} = 0$.

Note that while $F$ is quite smooth, being a primitive of an $H^3_x$ function, it does not necessarily decay at spatial infinity $|x| \to \infty$; in fact, it is not even bounded, since $H^3_x$ functions need not be absolutely integrable.  However, this problem is purely in the low frequencies of $F$; it is clear that $P_{hi} F$ is bounded (and in $H^4_x$), since $\partial_x$ is essentially an isometry from $H^4_x$ to $H^3_x$ on these frequency ranges.  

The function $e^{-iF}$ is clearly bounded, but obviously does not decay at infinity.  However, its derivative $\partial_x e^{-iF} = -iF_x e^{-iF} = -\frac{i}{2} u e^{-iF}$ can easily be shown to be in $H^3_x$.  In particular we can meaningfully talk about Fourier projections of $e^{-iF}$ as long as we stay away from the frequency origin; for instance, $P_{hi} e^{-iF}$ is well defined.  We can also define projections such as $P_{lo} e^{-iF}$ by defining $P_{lo} e^{-iF} := e^{-iF} - P_{hi} e^{-iF}$.  Similar considerations apply of course to $F$ itself.

From \eqref{eq:F-u} it is clear that in order to get bounds on $u$ it will suffice to obtain appropriate bounds
on $F$ which are of one higher degree of regularity.  Inspired by the Cole-Hopf transformation mentioned earlier, we introduce the complex-valued field $w: \R \times \R \to \C$ defined by
\begin{equation}\label{w-def}
 w := P_{+hi}( e^{-iF} );
\end{equation}
note that this projection is well defined since $e^{-iF}$ has a derivative in $H^3_x$.  If we then expand out $w_t - i w_{xx}$ and apply \eqref{eq:F} and \eqref{iH}, we obtain 
\begin{equation}\label{w-eq}
\begin{split}
w_t + H w_{xx} &= w_t - i w_{xx} \\
 &= P_{+hi}( (-iF_t - F_{xx} + i F_x^2) e^{-iF} )\\
 & = P_{+hi}( (iH F_{xx} - F_{xx}) e^{-iF})\\
 & = -2 P_{+hi}( P_-(F_{xx}) e^{-iF}) \\
 & = -2 P_{+hi}( P_-(F_{xx}) w) - 2 P_{+hi}( P_-(F_{xx}) P_{lo}(e^{-iF}) );
\end{split}
\end{equation}
the last identity follows from \eqref{w-def} and the fact that $P_{+hi}( P_-(F_{xx}) P_{-hi}(e^{-iF}) )$ vanishes.
Of the two terms, the second term is much better than the first, indeed it, and both of its factors, must necessarily be very low frequency.  The more interesting term is the first one.

This gives us an equation for $w$ in terms of $F$.  Note that the very worst type of interaction in the Benjamin-Ono equation, in which a very low frequency is interacted with a very high frequency, is almost completely absent in this equation, since both derivatives will fall on the low frequency.  (The $F_{xx}$ term must have lower frequency than $w$, otherwise $P_-(F_{xx}) w$ will have negative frequency and thus vanish when $P_{+hi}$ is applied.)

Now that we have a good equation for $w$, we must go back and obtain a good equation for $F$.  Define $F_{LO} := P_{LO} F$, $F_{+HI} := P_{+HI} F$, and $F_{-HI} := P_{-HI} F$, thus $F = F_{LO} + F_{+HI} + F_{-HI}$.  Since $F$ is real, we observe that $F_{LO}$ is also real, and that $F_{-HI} = \overline{F_{+HI}}$.  We differentiate \eqref{w-def}
and compute
\begin{align*}
w_x &= -i P_{+hi}( F_x e^{-iF} )\\
 &= -i P_{+hi}( (F_{+HI})_x e^{-iF} ) - i P_{+hi}( (F_{LO} + F_{-HI})_x e^{-iF} ) \\
 &= -i (F_{+HI})_x e^{-iF} + i (P_{lo}+P_{-hi})( (F_{+HI})_x e^{-iF} ) - i P_{+hi}( (F_{LO} + F_{-HI})_x e^{-iF} ) 
\end{align*}
and hence, multiplying by $ie^{iF}$ and rearranging, we obtain an equation
\begin{equation}\label{Fhi-eq}
 \partial_x F_{+HI} = i e^{iF} w_x + e^{iF} E
\end{equation}
for $F_{+HI}$, where the error term $E$ is given explicitly by
\begin{equation}\label{E-eq}
 E :=  (P_{lo} + P_{-hi})( e^{-iF} (F_{+HI})_x ) - 
P_{+hi}( (F_{-HI})_x e^{-iF})
P_{+hi}( (F_{LO})_x e^{-iF}).
\end{equation}
One can think of the formula \eqref{Fhi-eq} as a version of the Bony linearization formula for $e^{iF}$.

To prove \eqref{u-small} it will suffice by \eqref{eq:F-u} to prove the estimates
\begin{equation}\label{after-bootstrap}
\partial_x P_\pm F = O_{S^1}(\eps^2); \quad w = O_{S^2}(\eps^2).
\end{equation}
By a standard continuity argument using the $H^3_x$ global well-posedness theory, it will suffice to prove this estimate assuming the bootstrap hypothesis
\begin{equation}\label{before-bootstrap}
\partial_x P_\pm F = O_{S^1}(\eps); \quad w = O_{S^2}(\eps).
\end{equation}
We now deduce \eqref{after-bootstrap} from \eqref{before-bootstrap}.  We first deal with the low frequencies $F_{LO}$, which are easy to handle.  From the $L^2$ conservation law \eqref{l2-conserv} and the smallness of the initial data in $H^1_x$ (and hence in $L^2_x$) we clearly have $ u = O_{L^\infty_t L^2_x}(\eps^2)$;
applying $P_{LO}$ and then Bernstein's inequality (and H\"older in time) we thus see that
\begin{equation}\label{flo-est}
 \partial_x P_\pm F_{LO} = O_{S^1}(\eps^2)
\end{equation}
as desired.

Now we verify \eqref{after-bootstrap} for $w$.  We begin with the initial data $w(0)$.  From \eqref{eq:F-u} and the hypothesis on $u$, we have
$ F_x(0) = O_{H^1_x}(\eps^2)$,
and so from the chain rule and Sobolev (and the fact that $e^{-iF(0)}$ is bounded, since $F(0)$ is real) 
it is easy to see that
$\partial_x e^{-iF(0)} = O_{H^1_x}(\eps^2)$.
From \eqref{w-def} we thus have
$$ w(0) = O_{H^2_x}(\eps^2).$$
In order to verify \eqref{after-bootstrap} for $w$, it thus suffices by \eqref{strichartz} to verify that
$$ w_t + H w_{xx} = O_{L^1_t H^2_x}(\eps^2).$$
By \eqref{w-eq} it thus suffices to show that
$$ P_{+hi}(P_-(F_{xx}) w), P_{+hi}(P_-(F_{xx}) P_{lo}(e^{-iF}) ) = O_{L^1_t H^2_x}(\eps^2).$$

We first consider the term $P_{+hi}(P_-(F_{xx}) w)$.  From \eqref{before-bootstrap} we see that $P_-(F_{xx}) = O_{L^4_t L^\infty_x}(\eps)$ and $w = O_{L^\infty_t H^2_x}(\eps)$.  So the claim will follow from a H\"older in time, and the following paraproduct estimate.

\begin{lemma}\label{para}  For any functions $f(x)$ and $g(x)$, we have 
$$ \| P_{+hi}( P_-(f) g ) \|_{H^2_x} \lesssim \| P_-(f) \|_{L^\infty_x} \| g \|_{H^2_x}.$$
\end{lemma}

\begin{proof}  For any integer $k \geq 0$, observe that
$$ P_+ P_k( P_-(f) g ) = P_+ P_{k}( P_-(f) P_{\geq k-1} g ) $$
and hence (by the $L^2$ boundedness of $P_+ P_{k}$ and H\"older)
$$ \| P_+ P_k( P_-(f) g ) \|_{L^2_x} \lesssim \| P_-(f) \|_{L^\infty_x} \| P_{\geq k-1} g \|_{L^2_x}.$$
Multiplying by $2^{2k}$, and then square summing in $k$, we obtain
$$ \| P_{+hi}( P_-(f) g ) \|_{H^2_x} \lesssim \| P_-(f) \|_{L^\infty_x} (\sum_k (2^{2k} \| P_{\geq k-1} g \|_{L^2_x})^2)^{1/2}.$$
But by Plancherel we readily see that the latter factor on the right-hand side is $O(\| g \|_{H^2_x})$, and the claim follows.
\end{proof}

Now we estimate $P_{+hi}(P_-(F_{xx}) P_{lo}(e^{-iF}) )$.  Observe from all the frequency localizations that we may freely replace $F_{xx}$ by $P_{lo} F_{xx}$.  But from \eqref{flo-est} and frequency localization we see that $P_- P_{lo} F_{xx} = O_{L^\infty_t H^2_x}(\eps^2)$, while the boundedness of $e^{-iF}$ and frequency localization we have $P_{lo} e^{-iF} = O_{L^\infty_t C^2_x}(1)$.  Multiplying these two estimates we see that $P_{+hi}(P_-(F_{xx}) P_{lo}(e^{-iF}) ) = O_{L^1_t H^2_x}(\eps^2)$ as desired.

Finally, we verify \eqref{after-bootstrap} for $F_{HI}$. By conjugation invariance of the $S^k$ spaces it suffices to do this for $F_{+HI}$, since $F_{-HI} = \overline{F_{+HI}}$.  By \eqref{Fhi-eq} it thus suffices to show that
$$ e^{iF} w_x, e^{iF} E = O_{S^1}(\eps^2).$$
From \eqref{before-bootstrap} and Sobolev we know that $F_x = O_{L^\infty_t H^1_x}(\eps) = O_{L^\infty_t C^0_x}(\eps)$, and hence $e^{\pm iF} = O_{L^\infty_t C^1_x}(1)$.  In particular, we see that multiplication by $e^{\pm iF}$ is a bounded operation on $S^1$.  Thus it will suffice to show that
$$ w_x, E = O_{S^1}(\eps^2).$$
The claim for $w_x$ has already been proven, so it suffices to estimate $E$.  We expand $E$ using \eqref{E-eq}.  Consider the first term $(P_{lo} + P_{-hi})( e^{-iF} (F_{+HI})_x )$ in \eqref{E-eq}.  Note that we may freely replace $e^{-iF}$ by $P_{-hi} e^{-iF}$ due to all the frequency projections.  From \eqref{before-bootstrap} we have $F_x = O_{L^\infty_t H^1_x}(\eps)$, and hence $\partial_x e^{-iF} = -iF_x e^{-iF} = O_{L^\infty_t H^1_x}(\eps)$.  Hence $P_{-hi} e^{-iF} = O_{L^\infty_t H^2_x}(\eps) = O_{L^\infty_t C^1_x}(\eps)$ by Sobolev.  Since $(F_{+HI})_x = O_{S^1}(\eps)$, we thus easily see from Leibnitz and H\"older that
$$ (P_{lo} + P_{-hi})( e^{-iF} (F_{+HI})_x ) = O_{S^1}(\eps^2)$$
as desired (note that $P_{lo}$ and $P_{-hi}$ are bounded on $S^1$).  

The second term $P_{+hi}( (F_{-HI})_x e^{-iF})$ of \eqref{E-eq} is treated similarly, so we turn now to
$P_{+hi}( (F_{LO})_x e^{-iF} )$.  From \eqref{flo-est} we have that $(F_{LO})_x = O_{S^1}(\eps^2)$; since multiplication by $e^{-iF}$ is bounded on $S^1$, we see that $P_{+hi}( (F_{LO})_x e^{-iF} ) = O_{S^1}(\eps^2)$ as desired.
This completes the proof of Theorem \ref{apriori}.
\end{proof}

\section{The Lipschitz bound}\label{sec:lip}

The objective of this section is to show that the solution map is Lipschitz continuous from  $H^3_x \cap B(0,R)$, in the $L^2$ topology, to  $C^0_{[-T,T]} L^2_x$.  More precisely:

\begin{theorem}\label{lip-thm}  Let $0 < R \ll \eps^2$.  Then \eqref{lip} holds for all $u_0, \tilde u_0 \in H^3_x \cap B(0,R)$.
\end{theorem}

\begin{proof}  This will be an easy energy method argument.  Write $u(t) = S(t) u_0$ and $\tilde u(t) = S(t) \tilde u_0$, and let $v(t) := \tilde u(t) - u(t)$ denote the difference between the two solutions.  Then $v$ obeys the difference equation
$$
v_t + H v_{xx} = \frac{1}{2} ( \tilde u^2 - u^2 )_x = ( U v )_x
$$
where $U := \frac{\tilde u + u}{2}$.  Now we compute the growth of the $L^2$ norm of $v$:
\begin{align*}
\partial_t \int v(t)^2\ dx &= 2 \int v(t) v_t(t)\ dx \\
&= 2 \int v(t) (v_t(t) + H v_{xx}(t))\ dx \\
&= 2 \int v(t) (U(t)v(t))_x\ dx \\
&= - 2 \int v_x(t) U(t) v(t)\ dx \\
&= - \int (v^2(t))_x U(t)\ dx \\
&= \int v^2(t) U_x(t)\ dx \\
&= O(\| U_x(t) \|_{L^\infty_x}) \int v(t)^2\ dx.
\end{align*}
By Gronwall's inequality, we thus have
$$ \sup_{t \in [-1,1]} \int v(t)^2\ dx \leq \exp(\| U_x \|_{L^1_t L^\infty_x}) \int v(0)^2\ dx.$$
However, from Theorem \ref{apriori} we have
$ U_x = O_{L^4_t L^\infty_x}(\eps^2) = O_{L^1_t L^\infty_x}(\eps^2)$
and hence
$$ \sup_{t \in [-1,1]} \int v(t)^2\ dx \lesssim \int v(0)^2\ dx.$$
The claim \eqref{lip} then follows.
\end{proof}

From Theorem \ref{lip-thm} we know that the solution map $S(t)$ is a Lipschitz map from $H^3_x \cap B(0,R)$, with the $L^2_x$ topology, to $C^0_t H^3_x([-1,1] \times \R)$, with the $L^2_x$ topology.  Of course, in the $L^2_x$ topology, $H^3_x \cap B(0,R)$ is dense in $B(0,R)$.  Thus we may uniquely complete the solution map to be a Lipschitz continuous map from $B(0,R)$ with the $H^1_x$ topology, to $C^0_t H^1_x$ with the $H^1_x$ topology; this of course proves \eqref{lip} for all choices of initial data in $B(0,R)$.   Also we see from \eqref{lip} that the solution operator $S(t)$ constructed this way is unique in the class of limits of smooth solutions.

However, we are still missing several components of Theorem \ref{thm:lwp}.  Firstly, we have not yet shown that $S(t)$ is continuous in the $H^1_x$ topology (as opposed to the $L^2_x$ topology).  Also, we have not yet shown that the solution operator constructed by the above limiting procedure is a solution in the sense of distributions.  Finally, we have not shown the persistence of regularity bound \eqref{eq:regularity}.  To obtain all of these additional properties of the solution map $S(t)$, we will revisit the a priori bound in Theorem \ref{apriori}, but this time we will introduce the machinery of \emph{frequency envelopes}, as used recently in wave maps in \cite{tao:wavemap}, \cite{tao:wavemap2}, \cite{tataru}.

\section{Frequency envelopes}

In order to expedite the persistence of regularity and continuous dependence of the data results, we shall recall\footnote{Our definition of frequency envelope has been modified for our subcritical setting, as opposed to the critical setting of the papers cited above.  In particular, the frequency parameter $k$ is restricted to non-negative values.} the notion of a \emph{frequency envelope} (introduced in \cite{tao:wavemap}, \cite{tao:wavemap2} (see also \cite{kr}, \cite{mn}) and refined slightly in \cite{tataru} for applications to persistence of regularity).

Let $M$ be an integer larger than $s_0$, and let $0 < \sigma \ll 1$ be a small number depending on $M$ and $s_0$ to be chosen later.  We allow our implicit constants in the notation $X \lesssim Y$ to depend on $M$ and $\sigma$ as well as $s_0$.

\begin{definition}  A \emph{frequency envelope} $c$ is a map from positive integers $k$ to non-negative reals $c_k$ such that 
\begin{equation}\label{eq:c-norm}
c_0 \sim 1,
\end{equation}
\begin{equation}\label{eq:c-energy}
 \sum_{k=0}^\infty c_k^2 \lesssim 1
\end{equation}
and we have the log-Lipschitz conditions
\begin{equation}\label{eq:raise}
 2^{-Mr}c_k \leq c_{k+r} \leq 2^{\delta r} c_k
\end{equation}
for all $k,r \geq 0$.  
\end{definition}

From \eqref{eq:c-energy} we see in particular that the $c_k$ are bounded:
\begin{equation}\label{eq:cj-bound}
c_k \lesssim 1.
\end{equation}
The extra condition \eqref{eq:c-norm} is a technical one which simplifies the treatment of the low frequencies (in particular, it allows one to exploit the $L^2$ conservation law easily).  However, it is not strictly necessary and one can relax the notion of frequency envelope by dropping \eqref{eq:c-norm} (but then one has to estimate the low frequencies by iteration methods rather than via the conservation law).
The condition \eqref{eq:raise} allows the coefficients $c_k$ to drop quite fast as $k$ increases, but only allows them to increase very slowly.

We now apply this concept of a frequency envelope to translation-invariant Banach spaces.  If $X$ is a translation-invariant Banach space, either of functions of space $x$ or spacetime $(t,x)$, and $c$ is a frequency envelope, we define the envelope norm $X_c$ by
\begin{equation}\label{eq:hsc-def}
 \| f \|_{X_c} := \|f \|_X + \sup_{k \geq 0} \frac{\| P_k f \|_{X}}{c_k}.
\end{equation}
It is easy to verify that $X_c$ is indeed a norm.  One can think of $X_c$ as a Besov-type variant of $X$ which contains more refined information about the distribution of the ``energy'' of $f$ (as measured using the $X$ norm) in frequency space.

We now prove some linear and bilinear estimates involving envelope spaces that we will need in the next section.  First observe that the Strichartz estimate \eqref{strichartz} automatically implies the envelope variant
\begin{equation}\label{strichartz-envelope}
\| u \|_{S^k_c} \lesssim \| u(0) \|_{H^k_c} + \| u_t + H u_{xx} \|_{(L^1_t H^k_x)_c};
\end{equation}
this basically follows by applying \eqref{strichartz} to the Littlewood-Paley components $P_j u$ of $u$ and observing that the operator $\partial_t + H \partial_{xx}$ commutes with Littlewood-Paley operators.  In a similar vein we see that the spaces $S^k_c$ react well under differentiation, for instance we have
$$ \| u \|_{S^2_c} \sim \| u_x \|_{S^1_c} + \| u \|_{S^1_c},$$
and we can omit the lower-order term $\| u \|_{S^1_c}$ if $u$ is restricted to high frequencies $|\xi| \gtrsim 1$.

To prove some bilinear estimates, we shall need the following abstract lemma which allows us to convert bilinear estimates on Banach spaces $X,Y,Z$ to bilinear estimates on their envelope counterparts $X_c$, $Y_c$, $Z_c$.

Let $X,Y,Z$ be translation-invariant Banach spaces, and let $T: X \times Y \to Z$ be a bilinear operator.  We say that $T$ is of \emph{product type} if it respects the Littlewood-Paley trichotomy; in other words, that $P_{k} T(P_{k_1} f, P_{k_2} g)$ vanishes unless we are in one of the following three cases:

\begin{itemize}

\item (Low-High interactions) $k = k_2 + O(1)$ and $k_1 \leq k + O(1)$.

\item (High-Low interactions) $k = k_1 + O(1)$ and $k_2 \leq k + O(1)$.

\item (High-High interactions) $k_1 \geq k + O(1)$ and $k_2 = k_1 + O(1)$.

\end{itemize}

It is easy to see that the product operator $fg$ is of product type; it is clear that one can also add derivatives or frequency projections to the bilinear form $fg$ and still retain product type, for instance $P_{+hi}( P_-(f_x) g)$ is of product type.

\begin{lemma}\label{high} Let $X$, $Y$, $Z$ be translation-invariant Banach spaces, and let $c$ be a frequency envelope.  Suppose that we have a bilinear operator $T$ of product type which obeys the bilinear estimates
\begin{equation}\label{t}
 \| T(f,g) \|_Z \lesssim \|f\|_X \| g \|_Y
\end{equation}
for all test functions $f,g$.  Suppose also that we have the ``high-high refinement'' 
\begin{equation}\label{thh}
 \| P_{k} T( P_{k_1} f, P_{k_2} g ) \|_Z \lesssim 2^{-2\sigma (k_1-k)} \| P_{k_1} f \|_X \| P_{k_2} g \|_Y
\end{equation}
whenever $k_1 \geq k + O(1)$ and $k_2 = k_1 + O(1)$.  Then we have
$$ \| T(f,g) \|_{Z_c} \lesssim \|f\|_{X_c} \| g \|_{Y_c}.$$
\end{lemma}

\begin{proof}
Since
$$ \| T(f,g) \|_{Z} \lesssim \|f\|_X \|g\|_Y \lesssim \|f\|_{X_c} \| g \|_{Y_c}$$
it will suffice to prove the estimate
$$ \| P_k T(f,g) \|_Z \lesssim c_k \| f \|_{X_c} \| g \|_{Y_c}.$$
We split $f$ into Littlewood-Paley pieces $f = \sum_{k_1} P_{k_1} f$.  First consider the contribution where $k_1 = k + O(1)$.  Then by \eqref{t} we have
$$ \| P_k T(P_{k_1} f, g) \|_Z \lesssim \| T(P_{k_1} f, g) \|_Z \lesssim \| P_{k_1} f \|_{X} \| g \|_{Y}
\lesssim c_{k_1} \| f \|_{X_c} \| g \|_{Y_c}$$
which is acceptable since $c_{k_1} \sim c_k$ from \eqref{eq:raise}.  

Now consider the contributions where $k_1 < k - C$ for some large constant $C$.  Then the Littlewood-Paley trichotomy allows us to restrict $g$ to $\sum_{k_2 = k + O(1)} P_{k_2} g$, and we can estimate this contribution using the triangle inequality by
\begin{align*}
\sum_{k_2 = k + O(1)} \| P_k T(P_{<k-C} f, P_{k_2} g) \|_Z
&\lesssim \sup_{k_2 = k+O(1)} \| T(P_{<k-C} f, P_{k_2} g) \|_Z\\
&\lesssim \| P_{<k-C} f \|_X \sup_{k_2 = k+O(1)} \| P_{k_2} g \|_Y\\
&\lesssim \| f\|_{X_c} \sup_{k_2 = k + O(1)} c_{k_2} \| g \|_{Y_c}
\end{align*}
which is acceptable as before.

Finally, consider the contributions where $k_1 > k + C$ for some large constant $C$.  Then the Littlewood-Paley trichotomy allows us to restrict $g$ to $\sum_{k_2 = k_1 + O(1)} P_{k_2} g$, and we can estimate this contribution using the triangle inequality by
$$ \sum_{k_1 > k+C} \sum_{k_2 = k_1 + O(1)} \| P_k T(P_{k_1} f, P_{k_2} g) \|_Z.$$
Applying \eqref{thh} we can estimate this by
$$ \sum_{k_1 > k+C} \sum_{k_2 = k_1 + O(1)} 2^{-2\sigma(k_1-k)} \| P_{k_1} f \|_X \|P_{k_2} g \|_Y,$$
which we can estimate by
$$ \sum_{k_1 > k+C} \sum_{k_2 = k_1 + O(1)} 2^{-2\sigma(k_1-k)} c_{k_1} \| f \|_{X_c} \| g \|_{Y_c}.$$
But this is acceptable by \eqref{eq:raise}.
\end{proof}

Now we can state the bilinear estimates that we shall need.

\begin{proposition}\label{bilinear}  For any functions $f, g, \phi$ on $\R$, we have
\begin{equation}\label{algebra}
 \| f g \|_{H^j_c} \lesssim \|f \|_{H^2_c} \| g\|_{H^j_c}
\end{equation}
and
\begin{equation}\label{smooth}
 \| f \phi \|_{H^j_c} \lesssim \|f \|_{H^j_c} ( \| \phi \|_{L^\infty_x} + \| \phi_x \|_{H^{M+10}_x})
\end{equation}
for $k=1,2$.  We also have the spacetime variants
\begin{equation}\label{algebra-spacetime}
 \| f g \|_{(L^\infty_t H^j_x)_c} \lesssim \|f \|_{(L^\infty_t H^2_x)_c} \| g\|_{(L^\infty_t H^j_x)_c}
\end{equation}
and
\begin{equation}\label{smooth-spacetime}
 \| f \phi \|_{(L^\infty_t H^j_x)_c} \lesssim \|f \|_{(L^\infty_t H^j_x)_c} ( \| \phi \|_{L^\infty_t L^\infty_x} + \| \phi_x \|_{L^\infty_t H^{M+10}_x})
\end{equation}
where now $f,g,\phi$ are functions on $[-1,1] \times \R$.  We also have the variant
\begin{equation}\label{strichartz-multiplication}
 \| f \phi \|_{S^1_c} \lesssim \|f \|_{S^1_c} ( \| \phi \|_{L^\infty_t L^\infty_x} + \| \phi_x \|_{(L^\infty_t H^1_x)_c}).
\end{equation}
\end{proposition}

\begin{proof}
It suffices to prove the estimates \eqref{algebra-spacetime}, \eqref{smooth-spacetime}, \eqref{strichartz-multiplication}, since the claims \eqref{algebra}, \eqref{smooth} clearly follow as special cases.

We begin with \eqref{algebra-spacetime}.  Using Lemma \ref{high}, it suffices to prove the non-envelope space estimate
$$ \| f g \|_{L^\infty_t H^j_x} \lesssim \|f \|_{L^\infty_t H^2_x} \| g\|_{L^\infty_t H^j_x}$$
as well as the high-high refinement
$$ \| P_k( (P_{k_1} f) (P_{k_2} g) \|_{L^\infty_t H^j_x} \lesssim 2^{-2\sigma(k_1-k)} \|P_{k_1} f \|_{L^\infty_t H^2_x} \| P_{k_2} g\|_{L^\infty_t H^j_x}$$
whenever $k_1 = k_2 + O(1)$ and $k \leq k_1 + O(1)$.

The first inequality is clear from the Leibnitz, H\"older, and Sobolev inequalities, so we turn to the high-high refinement.  For this we use Bernstein's inequality to estimate
$$ \| P_k( (P_{k_1} f) (P_{k_2} g) \|_{L^\infty_t H^j_x} \lesssim 2^{jk} 2^{k/2} \| (P_{k_1} f) (P_{k_2} g) \|_{L^\infty_t L^1_x};$$
using H\"older we thus have
$$ \| P_k( (P_{k_1} f) (P_{k_2} g) \|_{L^\infty_t H^j_x} \lesssim 2^{jk} 2^{k/2} 2^{-2k_1} \| P_{k_1} f \|_{L^\infty_t H^2_x} 2^{-jk_2} \| P_{k_2} g \|_{L^\infty_t H^j_x}$$
and the claim readily follows (in fact we have plenty of powers of $2^k$ and $2^{k_1}$ to spare).

Now we prove \eqref{smooth-spacetime}.  A simple computation using Plancherel and \eqref{eq:c-norm}, \eqref{eq:raise} shows that
$$ \| \phi_x \|_{(L^\infty_t H^2_x)_c} \lesssim \| \phi_x \|_{L^\infty_t H^{M+10}_x}$$
so it suffices to prove the stronger estimate
\begin{equation}\label{sm-2}
 \| f \phi \|_{(L^\infty_t H^j_x)_c} \lesssim \|f \|_{(L^\infty_t H^j_x)_c} ( \| \phi \|_{L^\infty_t L^\infty_x} + \| \phi_x \|_{(L^\infty_t H^1_x)_c}).
\end{equation}
We split $\phi = P_{lo} \phi + P_{hi} \phi$.  Since
$\| P_{hi} \phi \|_{(L^\infty_t H^2_x)_c} \lesssim \| \phi_x \|_{(L^\infty_t H^1_x)_c}$
so we see from \eqref{algebra-spacetime} that the contribution of $P_{hi} \phi$ is acceptable.  Now consider the contribution of $P_{lo} \phi$.  We split $f P_{lo} \phi$ as $\sum_{k \geq 0} (P_k f) (P_{lo} \phi)$ and observe that the term $(P_k f) (P_{lo} \phi)$ has spatial Fourier support in the region $\{ \xi: \langle \xi \rangle \sim 2^k \}$.  Thus by \eqref{eq:raise} and dyadic decomposition it will suffice to prove the estimate
$$ \| (P_k f) (P_{lo} \phi) \|_{L^\infty_t H^j_x} \lesssim \| P_k f \|_{L^\infty_t H^j_x} \| \phi \|_{L^\infty_t L^\infty_x}$$
for each $k \geq 0$.  But we can replace the $H^j_x$ norm on both sides by an $L^2_x$ norm, losing a factor of $2^{jk}$ on both sides, and the claim then follows from H\"older (noting that $\| P_{lo} \phi \|_{L^\infty_t L^\infty_x} \lesssim \| \phi \|_{L^\infty_t L^\infty_x}$).  This proves \eqref{sm-2} and hence \eqref{smooth-spacetime}.

Finally, we prove \eqref{strichartz-multiplication}.  In light of \eqref{sm-2} it will suffice to show that
$$
 \| f \phi \|_{(L^4_t C^1_x)_c} \lesssim \|f \|_{(L^4_t C^1_x)_c} ( \| \phi \|_{L^\infty_t L^\infty_x} + \| \phi_x \|_{(L^\infty_t H^1_x)_c}).$$
Again we split $\phi = P_{lo} \phi + P_{hi} \phi$.  The treatment of $P_{lo} \phi$ is proven in exactly the same way as before.  To estimate $P_{hi} \phi$, it suffices by arguing as before to prove the estimate
$$ \| f g \|_{(L^4_t C^1_x)_c} \lesssim   \|f \|_{(L^4_t C^1_x)_c} \| g \|_{(L^\infty_t H^2_x)_c}.$$
We again apply Lemma \ref{high}.  The estimate
$$ \| f g \|_{L^4_t C^1_x} \lesssim   \|f \|_{L^4_t C^1_x} \| g \|_{L^\infty_t H^2_x}$$
can be easily deduced from Sobolev, H\"older, and the Leibnitz rule, so it suffices to prove the high-high refinement
$$ \| P_k( (P_{k_1} f) (P_{k_2} g)) \|_{L^4_t C^1_x} \lesssim  2^{-2\sigma(k_1 - k)}
\|P_{k_1} f \|_{L^4_t C^1_x} \| P_{k_2} g \|_{L^\infty_t H^2_x}$$
whenever $k_1 = k_2 + O(1)$ and $k \leq k_1 + O(1)$.  But this is clear, for instance we may estimate
\begin{align*} 
\|  P_k( (P_{k_1} f) (P_{k_2} g)) \|_{L^4_t C^1_x}
&\lesssim 2^k \| (P_{k_1} f) (P_{k_2} g)) \|_{L^4_t L^\infty_x} \\
&\lesssim 2^k \| P_{k_1} f \|_{L^4_t L^\infty_x} \| P_{k_2} g  \|_{L^\infty_t L^\infty_x} \\
&\lesssim 2^k 2^{-k_1} \| P_{k_1} f \|_{L^4_t C^1_x} \| P_{k_2} g  \|_{L^\infty_t H^2_x} 
\end{align*}
(in fact we can gain even more powers of $2^k$ and $2^{k_1}$ by being more careful).
\end{proof}

\section{A frequency-envelope a priori Strichartz bound}

We now prove a frequency envelope version of Theorem \ref{apriori}, namely:

\begin{theorem}\label{apriori-envelope}  Let $c$ be a frequency envelope, and let $u$ be an $H^{M+10}_x$ solution to \eqref{eq:benjamin} with $u(0) = O_{H^1_c}(\eps^2)$.  Then we have
$$ u = O_{S^1_c([-1,1] \times \R)}(\eps^2).$$
\end{theorem}

\begin{proof}  This will be a reprise\footnote{In fact, this theorem supercedes Theorem \ref{apriori}, which is actually redundant in this argument; but we elected to place the non-envelope version of the theorem first in order to clarify the exposition.} of the proof of Theorem \ref{apriori}.  Again, we restrict all spacetime norms to the slab $[-1,1] \times \R$.

Fix $c$, $u$.  We recall the fields $F, w, F_{LO}, F_{\pm HI}$ constructed in the proof of Theorem \ref{apriori}.
To prove the desired bound on $u$, it will suffice by \eqref{eq:F-u} to prove the estimates
\begin{equation}\label{after-bootstrap-envelope}
\partial_x P_\pm F = O_{S^1_c}(\eps^2); \quad w = O_{S^2_c}(\eps^2).
\end{equation}
By a standard continuity argument using the $H^{M+10}_x$ global well-posedness theory (noting that the $L^\infty_t H^{M+10}_x$ norm will control the $S^1_c$ and $S^2_c$ norms thanks to \eqref{eq:raise}), it will suffice to prove this estimate assuming the bootstrap hypothesis
\begin{equation}\label{before-bootstrap-envelope}
\partial_x P_\pm F = O_{S^1_c}(\eps); \quad w = O_{S^2_c}(\eps).
\end{equation}
We now deduce \eqref{after-bootstrap-envelope} from \eqref{before-bootstrap-envelope}.  From \eqref{flo-est}, \eqref{eq:c-norm}, and the low frequency projection we know that
\begin{equation}\label{flo-est-envelope}
 \partial_x F_{LO} = O_{S^1_c}(\eps^2)
\end{equation}
so the low frequencies are fine.

Now we verify \eqref{after-bootstrap-envelope} for $w$.  We begin with the initial data $w(0)$.  From \eqref{eq:F-u} we have
\begin{equation}\label{fx-envelope}
 F_x(0) = O_{H^1_c}(\eps^2)
\end{equation}
so in particular we have
$$ F_{HI}(0) = O_{H^2_c}(\eps^2).$$
From \eqref{algebra} and Taylor expansion we thus see (if $\eps$ is small) that
\begin{equation}\label{fhi}
 e^{-iF_{HI}(0)} = 1 + O_{H^2_c}(\eps^2).
\end{equation}
Now consider $F_{LO}(0)$.  From \eqref{fx-envelope} and frequency localization we have
$$ \partial_x F_{LO}(0) = O_{H^{M+10}_x}(\eps^2)$$
and hence by many applications of the chain rule (and using the boundedness of $e^{-iF_{LO}(0)}$)
$$ \partial_x e^{-iF_{LO}(0)} = O_{H^{M+10}_x}(\eps^2).$$
We can now multiply both sides of \eqref{fhi} by $e^{-iF_{LO}(0)}$ and use \eqref{smooth} to obtain
$$ e^{-iF(0)} = e^{-iF_{LO}(0)} + O_{H^2_c}(\eps^2).$$
Multiplying by $-iF_x(0)$ using \eqref{fx-envelope}, \eqref{algebra}, \eqref{smooth} we obtain
\begin{equation}\label{eifx}
 \partial_x e^{-iF(0)} = -i F_x(0) e^{-iF(0)} = O_{H^1_c}(\eps^2).
\end{equation}
From \eqref{w-def} we thus have
$$ w(0) = O_{H^2_c}(\eps^2).$$
In order to verify \eqref{after-bootstrap-envelope} for $w$, it thus suffices by \eqref{strichartz-envelope} to verify that
$$ w_t + H w_{xx} = O_{(L^1_t H^2_x)_c}(\eps^2).$$
By \eqref{w-eq} it thus suffices to show that
$$ P_{+hi}(P_-(F_{xx}) w), P_{+hi}(P_-(F_{xx}) P_{lo}(e^{-iF}) ) = O_{(L^1_t H^2_x)_c}(\eps^2).$$
The expression $P_{+hi}(P_-(F_{xx}) P_{lo}(e^{-iF}) )$ was already shown in the proof of Theorem \ref{apriori} to
be $O_{L^1_t H^2_x}(\eps^2)$, and is also of low frequency (supported on the region $|\xi| \leq 4$).  Hence it is acceptable by \eqref{eq:c-norm}.  Thus it will suffice to consider the term $P_{+hi}(P_-(F_{xx}) w)$.  From \eqref{before-bootstrap} we see that $P_-(F_{xx}) = O_{(L^4_t L^\infty_x)_c}(\eps)$ and $w = O_{(L^\infty_t H^2_x)_c}(\eps)$.  So the claim will follow from the following paraproduct estimate.

\begin{lemma}  For any functions $f(t,x)$ and $g(t,x)$, we have 
$$ \| P_{+hi}( P_-(f) g ) \|_{(L^1_t H^2_x)_c} \lesssim \| P_-(f) \|_{(L^4_t L^\infty_x)_c} \| g \|_{(L^\infty_t H^2_x)_c}.$$
\end{lemma}

\begin{proof}  We use Lemma \ref{high}.  From Lemma \ref{para} and a H\"older in time we already have
$$ \| P_{+hi}( P_-(f) g ) \|_{L^1_t H^2_x} \lesssim \| P_-(f) \|_{L^4_t L^\infty_x} \| g \|_{L^\infty_t H^2_x}.$$
It thus suffices to show the high-high refinement
$$ \| P_k P_{+hi}( (P_- P_{k_1} f) P_{k_2} g ) \|_{L^1_t H^2_x} \lesssim 2^{-2\sigma (k_1-k)} \| P_- P_{k_1}(f) \|_{L^4_t L^\infty_x} \| P_{k_2} g \|_{L^\infty_t H^2_x}$$
when $k_2 = k_1 + O(1)$ and $k \leq k_1 + O(1)$.  But we can compute
\begin{align*}
\| P_k P_{+hi}( (P_- P_{k_1} f) P_{k_2} g ) \|_{L^1_t H^2_x} &\lesssim 
2^{2k} \| (P_- P_{k_1} f) P_{k_2} g  \|_{L^1_t L^2_x} \\
&\lesssim 2^{2k} \| P_- P_{k_1} f \|_{L^4_t L^\infty_x} \| P_{k_2} g  \|_{L^\infty_t L^2_x} \\
&\lesssim 2^{2k} 2^{-2k_2} \| P_- P_{k_1} f \|_{L^4_t L^\infty_x} \| P_{k_2} g  \|_{L^\infty_t H^2_x} 
\end{align*}
and the claim follows.
\end{proof}

Finally, we verify \eqref{after-bootstrap-envelope} for $F_{HI}$. By conjugation invariance of the $S^k_c$ spaces it suffices to do this for $F_{+HI}$.  By \eqref{Fhi-eq} it thus suffices to show that
$$ e^{iF} w_x, e^{iF} E = O_{S^1}(\eps^2).$$
From \eqref{before-bootstrap} we know that $F_x = O_{(L^\infty_t H^1_x)_c}(\eps)$.  By repeating the analysis used to prove \eqref{eifx} (using \eqref{algebra-spacetime}, \eqref{smooth-spacetime} instead of \eqref{algebra}, \eqref{smooth}), we thus have
\begin{equation}\label{pmifx}
 \partial_x e^{\pm iF} = O_{(L^\infty_t H^1_x)_c}(\eps).
\end{equation}
By \eqref{strichartz-multiplication} we thus see that multiplication by $e^{\pm iF}$ is a bounded operation on $S^1_c$.
Thus it will suffice to show that
$$ w_x, E = O_{S^1_c}(\eps^2).$$
The claim for $w_x$ has already been proven, so it suffices to estimate $E$.  We expand $E$ using \eqref{E-eq}.  Consider the first term $(P_{lo} + P_{-hi})( e^{-iF} (F_{+HI})_x )$ in \eqref{E-eq}.  Note that we may freely replace $e^{-iF}$ by $P_{-hi} e^{-iF}$ due to all the frequency projections.  From \eqref{pmifx} we have
$$P_{-hi} e^{-iF} = O_{(L^\infty_t H^2_x)_c}(\eps).$$
In particular, from Sobolev we see that $P_{-hi} e^{-iF} = O_{L^\infty_t L^\infty_x}(\eps)$.  Since $(F_{+HI})_x = O_{S^1_c}(\eps)$, we thus see from \eqref{strichartz-multiplication} that
$$ (P_{lo} + P_{-hi})( e^{-iF} (F_{+HI})_x ) = O_{S^1_c}(\eps^2)$$
as desired.

The second term $P_{+hi}( (F_{-HI})_x e^{-iF})$ of \eqref{E-eq} is treated similarly, so we turn now to
$P_{+hi}( (F_{LO})_x e^{-iF} )$.  But this follows from \eqref{flo-est-envelope} since multiplication by $e^{-iF}$ is a bounded operation on $S^1_c$.  This completes the proof of Theorem \ref{apriori-envelope}.
\end{proof}

As a particular corollary of Theorem \ref{apriori-envelope}, we see that if $u$ is an $H^{M+10}_x$ solution of \eqref{eq:benjamin} with $u(0) = O_{H^1_c}(\eps^2)$, then we also have $u(t) = O_{H^1_c}(\eps^2)$ for all $t \in [-1,1]$, since the $S^1_c$ norm controls the $(L^\infty_t H^1_x)_c$ norm, which in turn controls the $L^\infty_t H^1_c$ norm.
This can be thought of a strengthened version of the energy conservation law \eqref{u-h1}; it asserts not only that the $H^1_x$ norm stays bounded, but in fact its frequency profile is also stable for short times.  

As an immediate consequence, we can prove \eqref{eq:regularity} for $H^{M+10}_x$ solutions.  Indeed, if $u$ is an $H^{M+10}_x$ solution, we can define $a_0, a_1, \ldots$ to be the sequence
$a_j := \| P_j u(0) \|_{H^1_x}$, and then let $c_j$ be the sequence
\begin{equation}\label{cj-def}
 c_j := 2^{-\sigma j} + \eps^{-2} ( \sum_{0 \leq k \leq j} 2^{-M|k-j|} a_k + \sum_{k > j} 2^{-\sigma|k-j|} a_k ).
\end{equation}
Then since $\| u(0) \|_{H^1_x} \lesssim \eps^2$, it is easy to verify from Plancherel that $c_j$ is a frequency envelope and
$$ \| u(0) \|_{H^1_c} \lesssim \eps^2.$$
Also, from \eqref{eq:planch} and Young's inequality we see that 
\begin{equation}\label{cj-l2}
 (\sum_j (2^{(s_0 - 1) j} c_j)^2)^{1/2} \lesssim  1 + \eps^{-2} (\sum_j (2^{(s_0 - 1) j} a_j)^2)^{1/2} \lesssim 1 + \eps^{-2} \| u(0) \|_{H^{s_0}_x}.
\end{equation}
Now from the corollary of Theorem \ref{apriori-envelope} mentioned earlier, we have
$$ \sup_{t \in [-1,1]} \| u(t) \|_{H^1_c} \lesssim \eps^2.$$
From this, \eqref{cj-l2}, and \eqref{eq:planch} we have
$$ \sup_{t \in [-1,1]} \| u(t) \|_{H^{s_0}_x} \lesssim \|u(0) \|_{H^{s_0}_x} + \eps^2$$
which is \eqref{eq:regularity} for $H^{M+10}_x$ solutions (note that by shrinking $\eps$ we can make
$\eps^2$ as small as $O(\|u(0) \|_{H^1_x})$, which is bounded by $O(\|u(0) \|_{H^{s_0}_x})$.  In the next section we shall establish continuity in $H^{s_0}_x$, which will allow us to remove this a priori $H^{M+10}_x$ restriction.

\section{Conclusion of proof of Theorem \ref{thm:lwp}}\label{conclude-sec}

Now we improve the $L^2_x$ continuity to $H^{s_0}_x$ continuity for any $s_0 \geq 1$, using the corollary to Theorem \ref{apriori-envelope}.  Let $u^{(n)}_0$ be any sequence of initial data in $H^{s_0}_x \cap B(0,R)$ which converges in $H^{s_0}_x$ norm to another initial datum $u_0$; we do not assume that $u_0$ is any smoother than $H^{s_0}_x$.  Once again, we may normalize $R \ll \eps^2$.  We already know from Theorem \ref{lip-thm} that $S(t) u^{(n)}_0$ converges to $S(t) u_0$ in the $C^0_t L^2_x$ topology; we wish to improve this to $C^0_t H^{s_0}_x$ convergence.  Our basic tool here is the frequency envelope spaces $H^{s_0}_c$, which are in some sense compact in $H^{s_0}_x$.

Let $0 < \alpha$.  We have to show that
\begin{equation}\label{eq:l2-conv}
 \sup_{t \in [-1,1]} \| S(t) u^{(n)}_0 - S(t) u_0 \|_{H^{s_0}_x} \lesssim \alpha
\end{equation}
for all sufficiently large $n$.  Since $\| u_0 \|_{H^{s_0}_x} < \infty$, we see from the monotone convergence theorem that there exists an integer $J \geq 1$ (depending\footnote{Note that this dependence on $u_0$ means that we establish merely that the solution map is continuous from $H^{s_0}_x$ to $H^{s_0}_x$; we know from the results of Koch-Tzvetkov \cite{kt-2} that the map cannot be \emph{uniformly} continuous in this topology.} on $u_0$, $\alpha$) such that
$$ \| P_{> J} u_0 \|_{H^{s_0}_x} \lesssim \alpha.$$
In particular, since $u^{(n)}_0$ converges to $u_0$ in $H^{s_0}_x$ we have
\begin{equation}\label{eq:j-lo}
 \| P_{> J} u^{(n)}_0 \|_{H^{s_0}_x} \lesssim \alpha
\end{equation}
for all sufficiently large $n$.  

Fix $n$ to be sufficiently large so that \eqref{eq:j-lo} holds.  Let $a_0, a_1, \ldots$ be the sequence (depending on $n$) defined by
$a_j := \| P_j u^{(n)}_0 \|_{H^1_x}$, and then let $c_j$ be the sequence defined by \eqref{cj-def}. 
Then since $\| u^{(n)}_0 \|_{H^1_x} \lesssim \eps^2$, it is easy to verify that $c_j$ is a frequency envelope and
$$ \| u^{(n)}_0 \|_{H^1_c} \lesssim \eps^2.$$
Also, from \eqref{eq:j-lo}, \eqref{eq:planch} we see that 
\begin{equation}\label{eq:c-hi}
 (\sum_{j \geq J+C} (2^{(s_0 - 1) j} c_j)^2)^{1/2} \lesssim \alpha/\eps^2,
\end{equation}
where the constant $C$ depends on $\sigma$, $\alpha$, $\eps$ but not on $n$.  By the corollary\footnote{Strictly speaking, this was only proved for $H^{M+10}_x$ solutions, but the claim follows for general solutions by a limiting argument using Fatou's lemma and the $L^2_x$ convergence result already proven.} to Theorem \ref{apriori} we thus have
$$ \sup_{t \in [-1,1]}\| S(t) u^{(n)}_0 \|_{H^1_c} \lesssim \eps^2.$$
In particular, by \eqref{eq:c-hi} we have
$$ \sup_{t \in [-1,1]}\| P_{>J}  S(t) u^{(n)}_0 \|_{H^{s_0}_x} \lesssim \alpha.$$
Taking limits as $n \to \infty$ and using Plancherel and Fatou's lemma, we obtain
$$ \sup_{t \in [-1,1]} \| P_{> J} S(t) u_0 \|_{H^{s_0}_x} \lesssim \alpha.$$
and hence
$$ \sup_{t \in [-1,1]}\| P_{> J} (S(t) u^{(n)}_0 - S(t) u_0) \|_{H^{s_0}_x} \lesssim \alpha.$$
Also, since $S(t) u^{(n)}_0 - S(t) u_0$ converges to zero in $C^0_t L^2_x$ norm from Theorem \ref{lip-thm}, we see that
$$ \sup_{t \in [-1,1]}\| P_{\leq J} (S(t) u^{(n)}_0 - S(t) u_0) \|_{H^{s_0}_x} \lesssim \alpha$$
for $n$ sufficiently large (depending on $\alpha$, $J$).  The claim \eqref{eq:l2-conv} then follows.

We have now shown that the solution map $S(t)$ can be uniquely continuously extended as a map from $H^{s_0} \cap B(0,R)$ (in the $H^{s_0}_x$ topology) to $H^{s_0}_x$; in particular we see by a standard limiting argument that the bound \eqref{eq:regularity}, which was already proven for $H^{M+10}_x \cap B(0,R)$ solutions, is in fact true for all $H^{s_0} \cap B(0,R)$ solutions.

Specializing to the case $s=1$, we have established that $S(t)$ is a continuous map from $B(0,R)$ to $H^1_x$ in the $H^1_x$ topology.  We can now easily take limits of \eqref{eq:benjamin} in the sense of distributions and conclude that \eqref{eq:benjamin} holds in the sense of distributions for all $u_0 \in B(0,R)$.  This proves all the conclusions of Theorem \ref{thm:lwp}, as desired.

\section{Remarks}

The methods employed here can certainly be extended below the $H^1_x$ regularity; we presented just the $H^1_x$ argument to simplify the exposition and to highlight the gauge transformation.  Notice that the only estimate we used for the linear equation $u_t + H u_{xx} = 0$ were the Strichartz estimates.  If instead one employed the Kato-type local smoothing estimate and the maximal function estimates as in \cite{kpv:kdv}, it is quite likely that one could lower the local existence result (or at least the a priori bound), probably to $H^{3/4+}_x$, in analogy with the situation for the KdV equation in \cite{kpv:kdv}.  The Lipschitz bound would then be more difficult to obtain, but perhaps one can use the parametrices for the linearized equation $v_t + H v_{xx} = ( U v )_x$ that were developed in \cite{kt}.  Alternatively one can attempt to gauge away the worst part of the non-linearity $(Uv)_x$ by techniques similar\footnote{One possibility is that one would have to microlocalize the gauge transformation, applying a different gauge to each Littlewood-Paley piece, in the spirit of \cite{tao:wavemap}, \cite{tao:wavemap2}, \cite{kr}, \cite{tataru}.} to that used for the a priori bound, e.g. by multiplying $v$ by $e^{-iF}$ where $F$ is some sort of primitive of $U$ as in \eqref{eq:F-u}.  Of course, as one lowers the regularity $H^{s_0}_x$ in which one obtains well-posedness, it may well be necessary to lower the regularity used for the Lipschitz bound as well.  Certainly one can only hope for a Lipschitz bound in a space no smoother than $H^{s_0-1/2}_x$, as the examples in \cite{kt-2} already show.

To go below $H^{3/4+}_x$, an obvious approach would be to introduce $X^{s,b}$ spaces, as used for instance in \cite{borg:xsb} or \cite{kpv}.  A preliminary analysis using the equation $w_t + i w_{xx} = P_{+hi}( P_-(\overline{w}_{xx}) w)$
as a toy model\footnote{Of course, if we wish to prove well-posedness for $u$ in $H^s_x$, then we expect to place $w$ in $H^{s+1}_x$, so for this discussion one would try to iterate the toy model equation for $w(0)$ in $H^{1+}_x$.} for \eqref{w-eq} suggests that such an approach should be able to obtain local well-posedness all the way down to $H^{0+}_x$, and perhaps even to $L^2_x$ if one uses some sort of Besov refinement of the $X^{s,b}$ spaces (e.g. the space $X^{0,1/2,1}$; see \cite{tataru:wave2}).  However, there is a non-trivial difficulty in employing $X^{s,b}$ spaces or related norms, because such norms do not react well to gauge transformations such as $v \mapsto e^{-iF} v$.  This is in contrast to Strichartz norms, which can easily absorb bounded factors such as $e^{-iF}$ with little difficulty.  It would be interesting to obtain a local well-posedness result at the $L^2_x$ regularity (which would immediately extend to a global well-posedness result, due to \eqref{l2-conserv}) but we do not know how to achieve this.  However it seems reasonably feasible to expect (for instance) a local well-posedness result at the level of $H^{1/2}_x$, which would then be global thanks to the conservation of the Hamiltonian \eqref{hamil}.  One may also be able to obtain global well-posedness at other regularities by use of such tools as the ``$I$-method'' (see e.g. \cite{ckstt:dnls}), although we do not pursue this question here.

The methods used here rely heavily on the fact that $u$ is real (in order for the gauge transformation $e^{-iF}$ to stay in $L^\infty$).  When $u$ is complex we suspect that the Benjamin-Ono equation is in fact ill-posed even in smooth topologies, unless one imposes some sort of moment condition on $u$ to eliminate the very low frequencies.

\end{document}